\documentclass[11pt]{article}

\usepackage{latexsym}
\usepackage{amssymb}

\newtheorem{thm}{Theorem}
\newtheorem{rmr}{Remark}

\begin{document}
{
\begin{center}
{\Large\bf
The Nevanlinna-type formula for the truncated matrix trigonometric moment problem.}
\end{center}
\begin{center}
{\bf S.M. Zagorodnyuk}
\end{center}

\section{Introduction.}
This paper is a continuation of our previous investigation on the truncated matrix trigonometric moment problem
by the operator approach in~\cite{cit_500_Z}.
The truncated matrix trigonometric moment problem
consists of finding a non-decreasing $\mathbb{C}_{N\times N}$-valued function
$M(t) = (m_{k,l})_{k,l=0}^{N-1}$, $t\in [0,2\pi]$, $M(0)=0$, which is
left-continuous in $(0,2\pi]$, and such that
\begin{equation}
\label{f1_1}
\int_0^{2\pi} e^{int} dM(t) = S_n,\qquad n=0,1,...,d,
\end{equation}
where $\{ S_n \}_{n=0}^d$ is a prescribed sequence of $(N\times N)$ complex matrices (moments).
Here $N\in \mathbb{N}$ and $d\in \mathbb{Z}_+$ are fixed numbers.
Set
\begin{equation}
\label{f1_2}
T_d = (S_{i-j})_{i,j=0}^d =
\left(
\begin{array}{ccccc} S_0 & S_{-1} & S_{-2} & \ldots & S_{-d}\\
S_1 & S_0 & S_{-1} & \ldots & S_{-d+1}\\
S_2 & S_1 & S_0 & \ldots & S_{-d+2}\\
\vdots & \vdots & \vdots & \ddots & \vdots\\
S_d & S_{d-1} & S_{d-2} & \ldots & S_0\end{array}
\right),
\end{equation}
where
$$ S_k := S_{-k}^*,\qquad k=-d,-d+1,...,-1, $$
and $\{ S_n \}_{n=0}^d$ are from~(\ref{f1_1}).
It is well known that the following condition:
\begin{equation}
\label{f1_4}
T_d\geq 0,
\end{equation}
is necessary and sufficient for the solvability of the moment problem~(\ref{f1_1}) (e.g.~\cite{cit_4500_A}).
The  moment problem~(\ref{f1_1}) is said to be {\it determinate} if it has a unique solution and {\it indeterminate}
in the opposite case.
We shall omit here an exposition on the history and recent results for the moment problem~(\ref{f1_1}).
All that can be found in~\cite{cit_500_Z}.

\noindent
The aim of our present investigation is to derive a Nevanlinna-type formula
for the truncated matrix trigonometric moment problem (TMTMP) in a general case.
Namely, we shall only assume that $d\geq 1$, condition~(\ref{f1_4}) is satisfied and the moment problem is
indeterminate, i.e. it has more than one solution.
The coefficients of the corresponding matrix linear  fractional transformation are
explicitly expressed by the prescribed moments.
Notice that in some situations (e.g. during a multiple application) the Nevanlina-type formula
has the advantage that the numbers of rows and the numbers of columns of its coefficients are less or equal to
$N$.
Easy conditions for the determinacy
of the moment problem in terms of the prescribed moments are given.

{\bf Notations. }
As usual, we denote by $\mathbb{R}, \mathbb{C}, \mathbb{N}, \mathbb{Z}, \mathbb{Z}_+$,
the sets of real numbers, complex numbers, positive integers, integers and non-negative integers,
respectively; $\mathbb{D} = \{ z\in \mathbb{C}:\ |z|<1 \}$.
The set of all complex matrices of size $(m\times n)$ we denote by $\mathbb{C}_{m\times n}$, $m,n\in \mathbb{N}$.
If $M\in \mathbb{C}_{m\times n}$ then $M^T$ denotes the transpose of $M$, and
$M^*$ denotes the complex conjugate of $M$. The identity matrix from $\mathbb{C}_{n\times n}$
we denote by $I_n$, $n\in \mathbb{N}$.

If H is a Hilbert space then $(\cdot,\cdot)_H$ and $\| \cdot \|_H$ mean
the scalar product and the norm in $H$, respectively.
Indices may be omitted in obvious cases.
By $\mathbb{C}^N$ we denote the finite-dimensional Hilbert space of complex column vectors
of size $N$ with the usual scalar product $(\vec x,\vec y)_{\mathbb{C}^N} = \sum_{j=0}^{N-1} x_j\overline{y_j}$,
for $\vec x,\vec y\in \mathbb{C}^N$,
$\vec x = (x_0,x_1,\ldots,x_{N-1})^T$, $\vec y = (y_0,y_1,\ldots,y_{N-1})^T$, $x_j,y_j\in \mathbb{C}$.

\noindent For a linear operator $A$ in $H$, we denote by $D(A)$
its  domain, by $R(A)$ its range, by $\mathop{\rm Ker}\nolimits A$
its null subspace (kernel), and $A^*$ means the adjoint operator
if it exists. If $A$ is invertible then $A^{-1}$ means its
inverse. $\overline{A}$ means the closure of the operator, if the
operator is closable. If $A$ is bounded then $\| A \|$ denotes its
norm.
For a set $M\subseteq H$
we denote by $\overline{M}$ the closure of $M$ in the norm of $H$.
If $M$ has a finite number of elements, then its number of elements we denote by
$\mathop{\rm card}\nolimits(M)$.
For an arbitrary set of elements $\{ x_n \}_{n\in I}$ in
$H$, we denote by $\mathop{\rm Lin}\nolimits\{ x_n \}_{n\in I}$
the set of all linear combinations of elements $x_n$,
and $\mathop{\rm span}\nolimits\{ x_n \}_{n\in I}
:= \overline{ \mathop{\rm Lin}\nolimits\{ x_n \}_{n\in I} }$.
Here $I$ is an arbitrary set of indices.
By $E_H$ we denote the identity operator in $H$, i.e. $E_H x = x$,
$x\in H$. In obvious cases we may omit the index $H$. If $H_1$ is a subspace of $H$, then $P_{H_1} =
P_{H_1}^{H}$ is an operator of the orthogonal projection on $H_1$
in $H$.

\section{The determinacy of the TMTMP. A Nevanlinna-type formula for the TMTMP.}
Let the moment problem~(\ref{f1_1}), with $d\geq 1$, be given and condition~(\ref{f1_4}),
with $T_d$ from~(\ref{f1_2}), be satisfied.
Let
$$ T_d = (\gamma_{n,m})_{n,m=0}^{(d+1)N-1},\
S_k = ( S_{k;s,l} )_{s,l=0}^{N-1},\quad -d\leq k\leq d, $$
where $\gamma_{n,m}, S_{k;s,l}\in \mathbb{C}$.
Observe that
\begin{equation}
\label{f2_9}
\gamma_{kN+s,rN+l} = S_{k-r;s,l},\qquad 0\leq k,r\leq d,\quad 0\leq s,l\leq N-1.
\end{equation}
We repeat here some constructions from~\cite{cit_500_Z}.
Consider a complex linear vector space $\mathfrak{H}$, which elements are row vectors
$\vec u = (u_0,u_1,u_2,...,u_{(d+1)N-1})$, with $u_n\in \mathbb{C}$, $0\leq n\leq (d+1)N-1$.
Addition and multiplication by a scalar are defined for vectors in a usual way.
Set
$$ \vec \varepsilon_n = (\delta_{n,0},\delta_{n,1},\delta_{n,2},...,\delta_{n,(d+1)N-1}),\qquad
0\leq n\leq (d+1)N-1, $$
where $\delta_{n,r}$ is Kronecker's delta.
In $\mathfrak{H}$ we define a linear functional $B$ by the following relation:
$$ B(\vec u, \vec w) = \sum_{n,r=0}^{(d+1)N-1} a_n \overline{b_r} \gamma_{n,r}, $$
where
$$ \vec u = \sum_{n=0}^{(d+1)N-1} a_n \vec\varepsilon_n,\quad
\vec w = \sum_{r=0}^{(d+1)N-1} b_r \vec\varepsilon_r,\quad a_n,b_r\in \mathbb{C}. $$
The space $\mathfrak{H}$ with $B$ form a quasi-Hilbert space (\cite{cit_6000_M}).
By the usual procedure of introducing of the classes of equivalence (see, e.g.~\cite{cit_6000_M}), we
put two elements $\vec u$, $\vec w$ from $\mathfrak{H}$ to the same class of equivalence denoted by
$[\vec u]$ or $[\vec w]$, if
$B(\vec u - \vec w,\vec u - \vec w) = 0$. The space of classes of equivalence is a (finite-dimensional) Hilbert
space. Everywhere in what follows it is denoted by $H$.
Set
$$ x_n := [\vec\varepsilon_n ],\qquad 0\leq n\leq (d+1)N-1. $$
Then
\begin{equation}
\label{f2_10}
(x_n,x_m)_H = \gamma_{n,m},\qquad 0\leq n,m\leq (d+1)N-1,
\end{equation}
and $\mathop{\rm span}\nolimits\{ x_n \}_{n=0}^{ (d+1)N-1 } = \mathop{\rm Lin}\nolimits\{ x_n \}_{n=0}^{ (d+1)N-1 }
= H$. Set $L_N := \mathop{\rm Lin}\nolimits\{ x_n \}_{n=0}^{ N-1 }$.
Consider the following operator:
\begin{equation}
\label{f2_11}
A x = \sum_{k=0}^{dN-1} \alpha_k x_{k+N},\quad x = \sum_{k=0}^{dN-1} \alpha_k x_k,\ \alpha_k\in \mathbb{C}.
\end{equation}
By~\cite[Theorem 3]{cit_500_Z} all solutions of the moment problem~(\ref{f1_1}) have the following form
\begin{equation}
\label{f2_42}
M(t) = (m_{k,j}(t))_{k,j=0}^{N-1},\qquad t\in [0,2\pi],
\end{equation}
where $m_{k,j}$ are obtained from the following relation:
\begin{equation}
\label{f2_43}
\int_0^{2\pi} \frac{1}{1-\zeta e^{it}} dm_{k,j} (t) =
(\left[
E_H - \zeta ( A \oplus \Phi_\zeta )
\right]^{-1} x_k,x_j)_H,\qquad \zeta\in \mathbb{D}.
\end{equation}
Here $\Phi_\zeta$ is an analytic in $\mathbb{D}$ operator-valued function which values are
linear contractions from $H\ominus D(A)$ into $H\ominus R(A)$.
Conversely, each analytic in $\mathbb{D}$ operator-valued function with above properties
generates by relations~(\ref{f2_42})-(\ref{f2_43}) a solution of the moment problem~(\ref{f1_1}).
The correspondence between all
analytic in $\mathbb{D}$ operator-valued functions with above properties
and all solutions of the moment problem~(\ref{f1_1}) is bijective.

Since we are going to obtain a Nevanlinna-type formula for the {\it indeterminate} TMTMP, it is important
to obtain some easy necessary and sufficient conditions for the determinacy of the TMTMP.

\begin{thm}
\label{t2_1}
Let the moment problem~(\ref{f1_1}), with $d\geq 1$, be given and condition~(\ref{f1_4}),
with $T_d$ from~(\ref{f1_2}), be satisfied. Let the operator $A$ in the Hilbert space $H$
be constructed as in~(\ref{f2_11}).
The following conditions are equivalent:

\begin{itemize}

\item[{\rm (A)}]
The moment problem~(\ref{f1_1}) is determinate;

\item[{\rm (B)}]
The defect numbers of $A$ are equal to zero;

\item[{\rm (C)}]
For each fixed number $r$, $dN\leq r\leq dN+N-1$, the following linear system of simultenuous equations:
\begin{equation}
\label{f2_44}
\sum_{n=0}^{dN-1} \alpha_{r,n} \gamma_{n,j} = \gamma_{r,j},
\end{equation}
with unknowns $\alpha_{r,0}, \alpha_{r,1}, \ldots, \alpha_{r,dN-1}$,
has a solution. Here the numbers $\gamma_{\cdot,\cdot}$ are defined by~(\ref{f2_9}).

\end{itemize}

If the above conditions are satisfied then the unique solution of the moment problem~(\ref{f1_1})
is given by the following relation:
\begin{equation}
\label{f2_45}
M(t) = (m_{k,j} (t))_{k,j=0}^{N-1},\quad
m_{k,j} (t) = (E_t x_k, x_j)_H,
\end{equation}
where $E_t$ is the left-continuous orthogonal resolution of unity of the unitary
operator $A$, which is piecewise constant.
\end{thm}
{\bf Proof. }
(A)$\Rightarrow$(B).
First, we notice that the defect numbers of $A$ are always equal, because $H$ is finite-dimensional and
$A$ is isometric. If the defect numbers are greater then zero, then we can choose unit vectors
$u_1\in H\ominus D(A)$ and $u_2\in H\ominus R(A)$. We set $\Phi_\zeta c u_1 = c u_2$, $\forall c\in \mathbb{C}$,
and
$\Phi_\zeta u = 0$, $u\in H\ominus D(A)$: $u\notin \mathop{\rm Lin}\nolimits\{ u_1 \}$;
$\zeta\in \mathbb{D}$. On the other hand, we set
$\widetilde \Phi_\zeta \equiv 0$. Functions $\Phi_\zeta$ and $\widetilde\Phi_\zeta$ produce
different solutions of the TMTMP by relation~(\ref{f2_43}).

\noindent
(B)$\Rightarrow$(A). If the defect numbers are zero, then the only admissible function $\Phi$
in relation~(\ref{f2_43}) is $\Phi_\zeta \equiv 0$.

Observe that (B) $\Leftrightarrow$ ($D(A)=H$) $\Leftrightarrow$
$$ \left(  x_{dN},x_{dN+1},\ldots,x_{dN+N-1} \in \mathop{\rm Lin}\nolimits\{ x_n \}_{n=0}^{dN-1} \right) $$
$\Leftrightarrow$
$$ \left( \left\{ \begin{array}{cccc}
x_{dN} = \sum_{n=0}^{dN-1} \alpha_{dN,n} x_n \\
x_{dN+1} = \sum_{n=0}^{dN-1} \alpha_{dN+1,n} x_n \\
\cdots \\
x_{dN+N-1} = \sum_{n=0}^{dN-1} \alpha_{dN+N-1,n} x_n \end{array}\right.,\
\alpha_{dN,n},\alpha_{dN+1,n},\cdots,\alpha_{dN+N-1,n}\in \mathbb{C} \right) $$
$\Leftrightarrow$
$$ \left( \left\{ \begin{array}{cccc}
(x_{dN},x_j)_H = (\sum_{n=0}^{dN-1} \alpha_{dN,n} x_n, x_j)_H \\
(x_{dN+1},x_j)_H = (\sum_{n=0}^{dN-1} \alpha_{dN+1,n} x_n, x_j)_H \\
\cdots \\
(x_{dN+N-1}, x_j)_H = (\sum_{n=0}^{dN-1} \alpha_{dN+N-1,n} x_n, x_j)_H \end{array}\right.\right., $$
where $\alpha_{dN,n},\alpha_{dN+1,n},\cdots,\alpha_{dN+N-1,n}\in \mathbb{C}$, $0\leq j\leq dN+N-1$)
$\Leftrightarrow$ (C).
$\Box$

We shall continue our considerations started before the statement of Theorem~\ref{t2_1}. In what follows
we assume that the TMTMP is indeterminate and the both defect numbers of $A$ are equal to $\delta = \delta(A)$,
$\delta\geq 1$.

\noindent
Let us apply the Gram-Schmidt orthogonalization procedure to the vectors $x_0,x_1,\ldots,x_{dN+N-1}$.
During this procedure we shall use the numbers $\gamma_{\cdot,\cdot}$ defined by~(\ref{f2_9}) and
the property~(\ref{f2_10}).

\noindent
{\bf Step $j$; $0\leq j\leq dN+N-1$.} Calculate
\begin{equation}
\label{f2_46}
n_{j} := \left\| x_{j} - \sum_{k:\ 0\leq k\leq j-1,\ n_k\not= 0} (x_j,y_k)_H y_k \right\|_H,
\end{equation}
where the sum on the right can be empty.
If $n_j\not= 0$, then we set
\begin{equation}
\label{f2_47}
y_{j} := \frac{1}{n_j} \left( x_{j} - \sum_{k:\ 0\leq k\leq j-1,\ n_k\not= 0} (x_j,y_k)_H y_k \right).
\end{equation}
If $n_j=0$, we pass to the next step.

\begin{rmr}
Notice that there exists a nonzero $n_j$ with $0\leq j\leq N-1$. In the opposite case
we would have $\| x_j \|_H^2 = \gamma_{j,j} = 0$,
$0\leq j\leq N-1$. Therefore $S_0=0$ and $M(t)\equiv 0$, and this contradicts to the indeterminacy of
the TMTMP.
\end{rmr}

\begin{rmr}
By~(\ref{f2_47}) every $y_j$ can be expressed as a linear combination of $x_0,x_1,\ldots,x_j$.
Thus, numbers $n_j$ can be calculated using the prescribed moments by relations~(\ref{f2_10}) and~(\ref{f2_9}).
\end{rmr}

Set $\Omega_1 = \{ j:\ 0\leq j\leq dN+N-1,\ n_j\not=0 \}$.
Then $\mathfrak{A} := \{ y_j \}_{j\in\Omega_1}$ is an orthonormal basis in $H$.
Moreover $\mathfrak{A}_1 := \{ y_j \}_{j\in\Omega_1:\ j\leq N-1}$ is an orthonormal basis in $L_N$, and
$\mathfrak{A}_2 :=\{ y_j \}_{j\in\Omega_1:\ j\leq dN-1}$ is an orthonormal basis in
$\mathop{\rm Lin}\nolimits\{ x_n \}_{n=0}^{dN-1} = D(A)$.
Therefore $\mathfrak{A}_3 :=\{ y_j \}_{j\in\Omega_1:\ dN\leq j\leq dN+N-1}$ is an orthonormal basis in
$H\ominus D(A)$. Consequently, $\delta\leq N$.

\noindent
Observe that $\mathop{\rm card}\nolimits(\mathfrak{A}_3) = \delta\geq 1$.
Set $\rho:= \mathop{\rm card}\nolimits(\mathfrak{A}_1)$,
$\tau:= \mathop{\rm card}\nolimits(\mathfrak{A}_2)$. Notice that $1\leq\rho\leq N$, $\tau\geq\rho\geq 1$.

\noindent
The $k$-th element, counting from zero, of the set $\mathfrak{A}$, arranged in the order of construction
of its elements, we denote by $u_k$, $k=0,1,...,\tau+\delta-1$.
Then $\mathfrak{A} = \{ u_k \}_{k=0}^{\tau+\delta-1}$,
$\mathfrak{A}_1 = \{ u_k \}_{k=0}^{\rho-1}$,
$\mathfrak{A}_2 = \{ u_k \}_{k=0}^{\tau-1}$,
$\mathfrak{A}_3 = \{ u_k \}_{k=\tau}^{\tau+\delta-1}$.

We need one more orthonormal basis in $H$.
Observe that
$\mathfrak{A}_2' := \{ v_k \}_{k=0}^{\tau-1}$, where $v_k := Au_k$, is an orthonormal basis in $R(A)$.
Notice that $R(A) = \mathop{\rm Lin}\nolimits\{ x_n \}_{n=dN}^{dN+N-1}$.
Therefore the linear span of vectors $\{ v_k \}_{k=0}^{\tau-1}$, $\{ x_n \}_{n=0}^{N-1}$,
is equal to $H$.
Hence, the linear span of vectors $\{ v_k \}_{k=0}^{\tau-1}$, $\{ u_k \}_{k=0}^{\rho-1}$,
is equal to $H$, as well.

\noindent
Let us apply the Gram-Schmidt orthogonalization procedure to the vectors
$v_0,v_1,\ldots, v_{\tau-1}, u_0,u_1,\ldots,u_{\rho-1}$.
Like in the above procedure, we shall use the numbers $\gamma_{\cdot,\cdot}$ defined by~(\ref{f2_9}) and
the property~(\ref{f2_10}).
Observe that the first $\tau$ elements are already orthonormal.

\noindent
{\bf Step $j$; $0\leq j\leq \rho-1$.} Calculate
\begin{equation}
\label{f2_48}
m_{j} := \left\| u_{j} - \sum_{l=0}^{\tau-1} (u_j,v_l)_H v_l -
\sum_{k:\ 0\leq k\leq j-1,\ m_k\not= 0} (u_j,f_k)_H f_k \right\|_H,
\end{equation}
where the last sum on the right can be empty.
If $m_j\not= 0$, then we set
\begin{equation}
\label{f2_49}
f_{j} := \frac{1}{m_j} \left( u_{j} - \sum_{l=0}^{\tau-1} (u_j,v_j)_H v_j -
\sum_{k:\ 0\leq k\leq j-1,\ m_k\not= 0} (u_j,f_k)_H f_k
\right).
\end{equation}
If $m_j=0$, we pass to the next step.

Set $\Omega_2 = \{ j:\ 0\leq j\leq \rho-1,\ m_j\not=0 \}$.
Then $\mathfrak{A}' := \{ v_k \}_{k=0}^{\tau-1} \cup \{ f_j \}_{j\in\Omega_2}$ is an orthonormal basis in $H$.
Set $\mathfrak{A}_3' := \{ f_j \}_{j\in\Omega_2}$.
Observe that $\mathop{\rm card}\nolimits(\mathfrak{A}_3') = \delta$.

\noindent
The $k$-th element, counting from zero, of the set $\mathfrak{A}_3'$, arranged in the order of construction
of its elements, we denote by $v_{\tau+k}$, $k=0,1,...,\delta-1$.
Then $\mathfrak{A}' = \{ v_k \}_{k=0}^{\tau+\delta-1}$,
$\mathfrak{A}_3' =\{ v_k \}_{k=\tau}^{\tau+\delta-1}$.

Denote by $\mathcal{M}_{1,\zeta}(\Phi)$ the matrix of the operator
$E_H - \zeta ( A \oplus \Phi_\zeta )$ in the basis $\mathfrak{A}$, $\zeta\in \mathbb{D}$. Here
$\Phi_\zeta$ is an analytic in $\mathbb{D}$ operator-valued function which values are
linear contractions from $H\ominus D(A)$ into $H\ominus R(A)$. Then
$$ \mathcal{M}_{1,\zeta}(\Phi) = \left( \left( \left[
E_H - \zeta ( A \oplus \Phi_\zeta ) \right] u_k, u_j
\right)_H
\right)_{j,k=0}^{\tau+\delta-1} =
\left(
\begin{array}{cc} A_{0,\zeta} & B_{0,\zeta}(\Phi) \\
C_{0,\zeta} & D_{0,\zeta}(\Phi) \end{array}
\right),
$$
where
$$ A_{0,\zeta} =
\left( \left( \left[
E_H - \zeta ( A \oplus \Phi_\zeta ) \right] u_k, u_j
\right)_H
\right)_{j,k=0}^{\tau-1}
=
\left( \left(
u_k - \zeta A u_k, u_j
\right)_H
\right)_{j,k=0}^{\tau-1} $$
\begin{equation}
\label{f2_50}
=
I_\tau - \zeta \left( \left(
v_k, u_j
\right)_H
\right)_{j,k=0}^{\tau-1},
\end{equation}
$$ B_{0,\zeta}(\Phi) =
\left( \left( \left[
E_H - \zeta ( A \oplus \Phi_\zeta ) \right] u_k, u_j
\right)_H
\right)_{0\leq j\leq \tau-1,\ \tau\leq k\leq \tau+\delta-1} $$
$$ =
\left( \left(
u_k - \zeta \Phi_\zeta u_k, u_j
\right)_H
\right)_{0\leq j\leq \tau-1,\ \tau\leq k\leq \tau+\delta-1} $$
$$ =
-\zeta \left( \left(
\Phi_\zeta u_k, u_j
\right)_H
\right)_{0\leq j\leq \tau-1,\ \tau\leq k\leq \tau+\delta-1}, $$
$$ C_{0,\zeta} =
\left( \left( \left[
E_H - \zeta ( A \oplus \Phi_\zeta ) \right] u_k, u_j
\right)_H
\right)_{\tau\leq j\leq \tau+\delta-1,\ 0\leq k\leq \tau-1} $$
$$ =
\left( \left(
u_k - \zeta A u_k, u_j
\right)_H
\right)_{\tau\leq j\leq \tau+\delta-1,\ 0\leq k\leq \tau-1} $$
\begin{equation}
\label{f2_51}
=
- \zeta \left( \left(
v_k, u_j
\right)_H
\right)_{\tau\leq j\leq \tau+\delta-1,\ 0\leq k\leq \tau-1},
\end{equation}
$$ D_{0,\zeta}(\Phi) =
\left( \left( \left[
E_H - \zeta ( A \oplus \Phi_\zeta ) \right] u_k, u_j
\right)_H
\right)_{\tau\leq j\leq \tau+\delta-1,\ \tau\leq k\leq \tau+\delta-1} $$
$$ =
\left( \left(
u_k - \zeta \Phi_\zeta u_k, u_j
\right)_H
\right)_{\tau\leq j\leq \tau+\delta-1,\ \tau\leq k\leq \tau+\delta-1} $$
$$ = I_\delta - \zeta
\left( \left(
\Phi_\zeta u_k, u_j
\right)_H
\right)_{\tau\leq j\leq \tau+\delta-1,\ \tau\leq k\leq \tau+\delta-1},\ \zeta\in \mathbb{D}. $$
Observe that the matrix $A_{0,\zeta}$ is invertible, since it is the matrix of the operator
$P_{D(A)} (E_H - \zeta A) P_{D(A)} = E_{D(A)} - \zeta P_{D(A)} A P_{D(A)}$, considered in
the Hilbert space $D(A)$, with respect to $\mathfrak{A}_2$, $\zeta\in \mathbb{D}$.
Notice that matrices $A_{0,\zeta},C_{0,\zeta}$, $\zeta\in \mathbb{D}$, can be calculated explicitly using
relations~(\ref{f2_10}) and~(\ref{f2_9}).

Denote by $F_\zeta$, $\zeta\in \mathbb{D}$, the matrix of the operator $\Phi_\zeta$,
acting from $H\ominus D(A)$ into $H\ominus R(A)$, with respect to the bases
$\mathfrak{A}_3$ and $\mathfrak{A}_3'$:
$$ F_\zeta = (f_\zeta(j,k))_{j,k=\tau}^{\tau+\delta-1},\qquad  f_\zeta(j,k) := (\Phi_\zeta u_k, v_j)_H. $$
Then
$$ \Phi_\zeta u_k = \sum_{l=\tau}^{\tau+\delta-1} f_\zeta(l,k) v_l,\quad \tau\leq k\leq \tau+\delta-1, $$
and
$$ B_{0,\zeta}(\Phi) =
-\zeta \left( \left(
\sum_{l=\tau}^{\tau+\delta-1} f_\zeta(l,k) v_l, u_j
\right)_H
\right)_{0\leq j\leq \tau-1,\ \tau\leq k\leq \tau+\delta-1} $$
$$ = -\zeta \left( \sum_{l=\tau}^{\tau+\delta-1} \left(
 v_l, u_j \right)_H f_\zeta(l,k)
\right)_{0\leq j\leq \tau-1,\ \tau\leq k\leq \tau+\delta-1},\quad \zeta\in \mathbb{D}. $$
Set
\begin{equation}
\label{f2_52}
W := \left( \left(  v_l, u_j \right)_H
\right)_{0\leq j\leq \tau-1,\ \tau\leq l\leq \tau+\delta-1}.
\end{equation}
Then
$$ B_{0,\zeta}(\Phi) = -\zeta W F_\zeta,\qquad \zeta\in \mathbb{D}. $$
We may write
$$ D_{0,\zeta}(\Phi) = I_\delta - \zeta
\left( \left(
\sum_{l=\tau}^{\tau+\delta-1} f_\zeta(l,k) v_l, u_j
\right)_H
\right)_{\tau\leq j\leq \tau+\delta-1,\ \tau\leq k\leq \tau+\delta-1} $$
$$ = I_\delta - \zeta
\left( \sum_{l=\tau}^{\tau+\delta-1}
\left( v_l, u_j \right)_H f_\zeta(l,k)
\right)_{\tau\leq j\leq \tau+\delta-1,\ \tau\leq k\leq \tau+\delta-1},\quad \zeta\in \mathbb{D}. $$
Set
\begin{equation}
\label{f2_53}
T :=
\left( \left( v_l, u_j \right)_H
\right)_{\tau\leq j\leq \tau+\delta-1,\ \tau\leq l\leq \tau+\delta-1}.
\end{equation}
Then
$$ D_{0,\zeta}(\Phi) = I_\delta - \zeta T F_\zeta,\qquad \zeta\in \mathbb{D}. $$
Thus, we may write
$$ \mathcal{M}_{1,\zeta}(\Phi) =
\left(
\begin{array}{cc} A_{0,\zeta} & -\zeta W F_\zeta \\
C_{0,\zeta} & I_\delta - \zeta T F_\zeta \end{array}
\right),\quad \zeta\in \mathbb{D},
$$
where $A_{0,\zeta}$, $C_{0,\zeta}$ are given by~(\ref{f2_50}),(\ref{f2_51}),
and $W,T$ are given by~(\ref{f2_52}),(\ref{f2_53}).
Let apply the Frobenius formula for the inverse of a block matrix~\cite[p. 59]{cit_7000_G}. Then
$$ \mathcal{M}_{1,\zeta}^{-1}(\Phi) =
\left(
\begin{array}{cc}  A_{0,\zeta}^{-1} -\zeta A_{0,\zeta}^{-1} W F_\zeta H_{\zeta}^{-1}(\Phi)
C_{0,\zeta} A_{0,\zeta}^{-1} & \ast \\
\ast & \ast \end{array}
\right), $$
$$ =
\left(
\begin{array}{cc}  \frac{1}{h_\zeta} A_{0,\zeta}^{+} - \frac{\zeta}{h_\zeta^2}
A_{0,\zeta}^{+} W F_\zeta H_{\zeta}^{-1}(\Phi)
C_{0,\zeta} A_{0,\zeta}^{+} & \ast \\
\ast & \ast \end{array}
\right),\quad \zeta\in \mathbb{D}, $$
where by stars ($\ast$) we denoted the blocks which are not of interest for us, and
$$ H_{\zeta}(\Phi) = I_\delta - \zeta T F_\zeta + \zeta C_{0,\zeta} A_{0,\zeta}^{-1} W F_\zeta
= I_\delta - \zeta T F_\zeta + \frac{\zeta}{h_\zeta} C_{0,\zeta} A_{0,\zeta}^{+} W F_\zeta $$
\begin{equation}
\label{f2_54}
= I_\delta +
\left(\frac{\zeta}{h_\zeta} C_{0,\zeta} A_{0,\zeta}^{+} W - \zeta T \right) F_\zeta,\quad \zeta\in \mathbb{D}.
\end{equation}
Here $A_{0,\zeta}^{+}$
denotes the adjoint matrix of $A_{0,\zeta}$, i.e. the transpose of the cofactor matrix, and
\begin{equation}
\label{f2_54_1}
h_\zeta = \det A_{0,\zeta}.
\end{equation}

\noindent
Let $\zeta\in \mathbb{D}$.
The minor of $\mathcal{M}_{1,\zeta}^{-1}(\Phi)$, standing in the first $\rho$ rows and the first $\rho$ columns,
we denote by $\mathcal{M}_{2,\zeta}(\Phi)$.
The minor of $A_{0,\zeta}^{+}$, standing in the first $\rho$ rows and the first $\rho$ columns,
we denote by $A_{1,\zeta}$. The first $\rho$ rows of $A_{0,\zeta}^{+}$ we denote by $A_{2,\zeta}$.
The first $\rho$ columns of $A_{0,\zeta}^{+}$ we denote by $A_{3,\zeta}$. Then
\begin{equation}
\label{f2_55}
\mathcal{M}_{2,\zeta}(\Phi) =
\frac{1}{h_\zeta} A_{1,\zeta} - \frac{\zeta}{h_\zeta^2}
A_{2,\zeta} W F_\zeta H_{\zeta}^{-1}(\Phi)
C_{0,\zeta} A_{3,\zeta},\quad \zeta\in \mathbb{D}.
\end{equation}
Observe that $\mathcal{M}_{2,\zeta}(\Phi)$ is the matrix of the operator
$P_{L_N} \left[ E_H - \zeta ( A \oplus \Phi_\zeta )
\right]^{-1} P_{L_N}$, considered as an operator in $L_N$, with respect to
the basis $\mathfrak{A}_1$, $\zeta\in \mathbb{D}$.

Consider the following operator from $\mathbb{C}^N$ to $L_N$:
$$ K \sum_{n=0}^{N-1} c_n \vec e_n = \sum_{n=0}^{N-1} c_n x_n,\qquad c_n\in \mathbb{C}, $$
where $\vec e_n = (\delta_{n,0},\delta_{n,1},\ldots,\delta_{n,N-1})\in \mathbb{C}^N$.
Let $\mathcal{K}$ be the matrix of $K$ with respect to the orthonormal bases
$\{ \vec e_n \}_{n=0}^{N-1}$ and $\mathfrak{A}_1$:
\begin{equation}
\label{f2_55_1}
\mathcal{K} = \left( \left(K \vec e_k, u_j \right)_H \right)_{0\leq j\leq \rho-1,\ 0\leq k\leq N-1}
= \left( \left(x_k, u_j \right)_H \right)_{0\leq j\leq \rho-1,\ 0\leq k\leq N-1}.
\end{equation}
Then we may write
$$ \left(\left[
E_H - \zeta ( A \oplus \Phi_\zeta )
\right]^{-1} x_k,x_j \right)_H $$
$$ =
\left( P_{L_N} \left[
E_H - \zeta ( A \oplus \Phi_\zeta )
\right]^{-1} P_{L_N} K \vec e_k, K \vec e_j \right)_{H} $$
$$ = \left( K^* P_{L_N} \left[
E_H - \zeta ( A \oplus \Phi_\zeta )
\right]^{-1} P_{L_N} K \vec e_k, \vec e_j \right)_{\mathbb{C}^N},\quad \zeta\in \mathbb{D}. $$
Observe that the right-hand side is equal to the element of the matrix
$\mathcal{K}^* \mathcal{M}_{2,\zeta}(\Phi) \mathcal{K}$, standing in row $j$, column $k$.
By~(\ref{f2_43}) we may write:
\begin{equation}
\label{f2_56}
\int_0^{2\pi} \frac{1}{1-\zeta e^{it}} dM^T (t) =
\mathcal{K}^* \mathcal{M}_{2,\zeta}(\Phi) \mathcal{K},\qquad \zeta\in \mathbb{D}.
\end{equation}
Set
$$
\mathbf{C}_{\zeta} = \zeta C_{0,\zeta} A_{0,\zeta}^{+} W - \zeta h_\zeta T,\quad
\mathbf{A}_{\zeta} = \mathcal{K}^* A_{1,\zeta} \mathcal{K}, $$
\begin{equation}
\label{f2_57}
\mathbf{B}_{\zeta} = \mathcal{K}^* A_{2,\zeta} W,\quad
\mathbf{D}_{\zeta} = C_{0,\zeta} A_{3,\zeta} \mathcal{K},\quad \zeta\in \mathbb{D}.
\end{equation}
By~(\ref{f2_54}),(\ref{f2_55}),(\ref{f2_56}) we get
$$
\int_0^{2\pi} \frac{1}{1-\zeta e^{it}} dM^T (t) =
\frac{1}{h_\zeta} \mathbf{A}_{\zeta} - \frac{\zeta}{h_\zeta^2}
\mathbf{B}_{\zeta} F_\zeta
\left(
I_\delta +
\frac{1}{h_\zeta} \mathbf{C}_{\zeta} F_\zeta
\right)^{-1}
\mathbf{D}_{\zeta},
$$
where $\zeta\in \mathbb{D}$.

\begin{thm}
\label{t2_2}
Let the moment problem~(\ref{f1_1}), with $d\geq 1$, be given and condition~(\ref{f1_4}),
with $T_d$ from~(\ref{f1_2}), be satisfied. Suppose that the moment problem is
indeterminate.
All solutions of the moment problem~(\ref{f1_1}) can be obtained from the following relation:
$$ \int_0^{2\pi} \frac{1}{1-\zeta e^{it}} dM^T (t) $$
\begin{equation}
\label{f2_59}
=
\frac{1}{h_\zeta} \mathbf{A}_{\zeta} - \frac{\zeta}{h_\zeta^2}
\mathbf{B}_{\zeta} F_\zeta
\left(
I_\delta +
\frac{1}{h_\zeta} \mathbf{C}_{\zeta} F_\zeta
\right)^{-1}
\mathbf{D}_{\zeta},\quad \zeta\in \mathbb{D},
\end{equation}
where $\mathbf{A}_{\zeta}$, $\mathbf{B}_{\zeta}$, $\mathbf{C}_{\zeta}$, $\mathbf{D}_{\zeta}$,
are matrix polynomials defined by~(\ref{f2_57}), with values in $\mathbb{C}_{N\times N}$, $\mathbb{C}_{N\times \delta}$,
$\mathbb{C}_{\delta\times \delta}$, $\mathbb{C}_{\delta\times N}$, respectively.
The scalar polynomial $h_\zeta$, $\deg h_\zeta\leq \tau$, is given by~(\ref{f2_54_1}).
Here $F_\zeta$ is an analytic in $\mathbb{D}$, $\mathbb{C}_{\delta\times\delta}$-valued function which values are
such that $F_\zeta^* F_\zeta \leq 1$, $\forall\zeta\in \mathbb{D}$.
Conversely, each analytic in $\mathbb{D}$, $\mathbb{C}_{\delta\times\delta}$-valued function
such that $F_\zeta^* F_\zeta \leq 1$, $\forall\zeta\in \mathbb{D}$,
generates by relation~(\ref{f2_59}) a solution of the moment problem~(\ref{f1_1}).
The correspondence between all
analytic in $\mathbb{D}$, $\mathbb{C}_{\delta\times\delta}$-valued functions such that $F_\zeta^* F_\zeta \leq 1$,
$\forall\zeta\in \mathbb{D}$,
and all solutions of the moment problem~(\ref{f1_1}) is bijective.
\end{thm}
{\bf Proof. }
The proof follows from the preceeding considerations.
$\Box$

\noindent
{\bf Example 2.1. } Let $N=3$, $d=1$,
$S_0 = \left(
\begin{array}{ccc} 1 & 1 & 0\\
1 & 1 & 0\\
0 & 0 & 1\end{array}
\right)$,
$S_1 = \left(
\begin{array}{ccc} 1 & 1 & 0\\
1 & 1 & 0\\
0 & 0 & 0\end{array}
\right)$. Consider the TMTMP with moments $S_0,S_1$.
It is straightforward to check that condition~(\ref{f1_4}) holds and condition (C) of Theorem~\ref{t2_1}
fails. Thus, the TMTMP is solvable and indeterminate.
The matrix $T_1$ from~(\ref{f1_2}) has the following form:
$$ T_1 = (\gamma_{n,m})_{n,m=0}^5 =
\left(
\begin{array}{cccccc}
1 & 1 & 0 & 1 & 1 & 0\\
1 & 1 & 0 & 1 & 1 & 0\\
0 & 0 & 1 & 0 & 0 & 0\\
1 & 1 & 0 & 1 & 1 & 0\\
1 & 1 & 0 & 1 & 1 & 0\\
0 & 0 & 0 & 0 & 0 & 1\end{array}
\right). $$
Let $H$ be the Hilbert space described above, after formula~(\ref{f2_10}), and $\{ x_n \}_{n=0}^5$ be elements with
the property~(\ref{f2_10}).

\noindent
Let us apply the orthogonalization procedure~(\ref{f2_46}),(\ref{f2_47}) to the elements
$x_0,x_1,x_2,x_3,x_4,x_5$.

\noindent
{\bf Step 0}. Calculate $n_0 = \| x_0 \|_H = \sqrt{(x_0,x_0)_H} = \sqrt{\gamma_{0,0}} = 1\not=0$.
Then we set $y_0 = x_0$.

\noindent
{\bf Step 1}. Calculate
$$ n_1^2 = \| x_1 - (x_1,y_0)_H y_0 \|_H^2 = \| x_1 - (x_1,x_0)_H x_0 \|_H^2
= (x_1 - \gamma_{1,0} x_0, x_1 - \gamma_{1,0} x_0)_H $$
$$ = (x_1 - x_0, x_1 - x_0)_H = (x_1,x_1)_H - (x_1,x_0)_H - (x_0,x_1)_H + (x_0,x_0)_H $$
$$ = \gamma_{1,1} - \gamma_{1,0} - \gamma_{0,1} + \gamma_{0,0} = 0. $$
Therefore we pass to the next step.

\noindent
{\bf Step 2}. We calculate
$$ n_2^2 = \| x_2 - (x_2,y_0)_H y_0 \|_H^2 = \| x_2 - \gamma_{2,0} x_0 \|_H^2 =
(x_2, x_2)_H = \gamma_{2,2} = 1. $$
Set
$$ y_2 = x_2 - (x_2,y_0)_H y_0 = x_2 - \gamma_{2,0} y_0 = x_2. $$
In step 3 we obtain $n_3=0$, in step 4 we get $n_4=0$.
Finally, in step~5 we get $n_5=1$, and $y_5=x_5$.

\noindent
Set $\mathfrak{A} = \{ y_0,y_2,y_5 \}$. Let $u_0:=y_0=x_0$, $u_1:=y_2=x_2$, $u_2:=y_5=x_5$.
Then $\mathfrak{A} = \{ u_k \}_{k=0}^2$. Observe that in our case we have: $\rho=\tau = 2$, $\delta=1$.

\noindent
Set $v_0 := Au_0 = Ax_0 = x_3$, $v_1 := Au_1 = Ax_2 = x_5$.
Let us apply the orthogonalization procedure~(\ref{f2_48}),(\ref{f2_49}) to the elements
$v_0,v_1,u_0,u_1$.

\noindent
{\bf Step 0}. Calculate
$$ m_0^2 =
\| u_0 - (u_0,v_0)_H v_0 - (u_0,v_1)_H v_1 \|_H^2
= \| x_0 - (x_0,x_3)_H x_3 - (x_0,x_5)_H x_5 \|_H^2 = $$
$$ \| x_0 - \gamma_{0,3} x_3 - \gamma_{0,5} x_5 \|_H^2
= \| x_0 - x_3 \|_H^2 = (x_0 - x_3,x_0-x_3)_H $$
$$ = \gamma_{0,0} - \gamma_{0,3} - \gamma_{3,0} + \gamma_{3,3} = 0. $$
Then we pass to the next step.

\noindent
{\bf Step 1}. Calculate
$$ m_1^2 = \| u_1 - (u_1,v_0)_H v_0 - (u_1,v_1)_H v_1 \|_H^2 = $$
$$ \| x_2 - (x_2,x_3)_H x_3 - (x_2,x_5)_H x_5 \|_H^2
= \| x_2 - \gamma_{2,3} x_3 - \gamma_{2,5} x_5 \|_H^2 $$
$$ = \| x_2 \|_H^2 = (x_2,x_2)_H = \gamma_{2,2} = 1. $$
Set
$$ f_1 = u_1 - (u_1,v_0)_H v_0 - (u_1,v_1)_H v_1 $$
$$ = x_2 - (x_2,x_3)_H x_3 - (x_2,x_5)_H x_5 = x_2 - \gamma_{2,3} x_3 - \gamma_{2,5} x_5 = x_2. $$
Set $\mathfrak{A}' = \{ v_0,v_1,f_1 \}$. Let $v_2 = f_1 = x_2$. Then
$\mathfrak{A}' = \{ v_k \}_{k=0}^2$.

\noindent
By~(\ref{f2_52}),(\ref{f2_53}) we may write
$$ W = \left( \left(  v_l, u_j \right)_H
\right)_{0\leq j\leq 1,\ 2\leq l\leq 2} =
\left(
\begin{array}{cc}
(x_2,x_0)_H \\
(x_2,x_2)_H \end{array}
\right)
=
\left(
\begin{array}{cc}
\gamma_{2,0} \\
\gamma_{2,2}\end{array}
\right)
=
\left(
\begin{array}{cc}
0 \\
1\end{array}
\right); $$
$$ T =
\left( \left( v_l, u_j \right)_H
\right)_{2\leq j\leq 2,\ 2\leq l\leq 2} = (v_2,u_2)_H = (x_2, x_5)_H = \gamma_{2,5} = 0. $$
By~(\ref{f2_50}),(\ref{f2_51}) we calculate:
$$ A_{0,\zeta} =
I_2 - \zeta \left( \left(
v_k, u_j
\right)_H
\right)_{j,k=0}^{1}
= I_2 - \zeta
\left(
\begin{array}{cc} (x_3,x_0)_H & (x_5,x_0)_H \\
(x_3,x_2)_H & (x_5,x_2)_H \end{array}
\right) $$
$$ =
\left(
\begin{array}{cc} 1-\zeta & 0 \\
0 & 1 \end{array}
\right); $$
$$ C_{0,\zeta} =
- \zeta \left( \left(
v_k, u_j
\right)_H
\right)_{2\leq j\leq 2,\ 0\leq k\leq 1} = - \zeta (\gamma_{3,5},\gamma_{5,5}) = -\zeta (0,1). $$
Then
$h_\zeta = \det A_{0,\zeta} = 1-\zeta$,
$A_{0,\zeta}^+ = \left(
\begin{array}{cc} 1 & 0 \\
0 & 1-\zeta \end{array}
\right) = A_{1,\zeta}=A_{2,\zeta}=A_{3,\zeta}$.
By~(\ref{f2_55_1}) we may write
$$ \mathcal{K} = \left( \left(x_k, u_j \right)_H \right)_{0\leq j\leq 1,\ 0\leq k\leq 2} =
\left(\begin{array}{ccc} \gamma_{0,0} & \gamma_{1,0} & \gamma_{2,0} \\
\gamma_{0,2} & \gamma_{1,2} & \gamma_{2,2} \end{array}\right)
=
\left(\begin{array}{ccc} 1 & 1 & 0 \\
0 & 0 & 1 \end{array}\right).
$$
Using~(\ref{f2_57}) we calculate
$$ \mathbf{C}_{\zeta} = -\zeta^2(1-\zeta),\quad
\mathbf{A}_{\zeta} =
\left(\begin{array}{ccc} 1 & 1 & 0 \\
1 & 1 & 0 \\
0 & 0 & 1-\zeta \end{array}\right),\quad
\mathbf{B}_{\zeta} =
\left(\begin{array}{ccc} 0 \\
0 \\
1-\zeta \end{array}\right), $$
$$ \mathbf{D}_{\zeta} = -\zeta (0,0,1-\zeta). $$
Finally, by~(\ref{f2_59}) we obtain:
$$ \int_0^{2\pi} \frac{1}{1-\zeta e^{it}} dM^T (t)
= \left(\begin{array}{ccc} \frac{1}{1-\zeta} & \frac{1}{1-\zeta} & 0 \\
\frac{1}{1-\zeta} & \frac{1}{1-\zeta} & 0 \\
0 & 0 & 1 + \zeta^2 \frac{F_\zeta}{1-\zeta^2 F_\zeta}\end{array}\right),\ \zeta\in \mathbb{D}. $$
In particular, if we choose $F_\zeta\equiv 1$, then
$$ M(t) =
\left(\begin{array}{ccc} \widetilde m(t) & \widetilde m(t) & 0 \\
\widetilde m(t) & \widetilde m(t) & 0 \\
0 & 0 & \widehat m(t) \end{array}\right),\ t\in [0,2\pi], $$
where
$$ \widetilde m(t) =
\left\{\begin{array}{cc} 0, & \mbox{if }t=0 \\
1, & \mbox{if }t\in(0,2\pi]\end{array}\right.,\quad
\widehat m(t) =
\left\{\begin{array}{ccc} 0, & \mbox{if }t=0 \\
\frac{1}{2}, & \mbox{if }t\in(0,\pi] \\
1, & \mbox{if }t\in(\pi,2\pi]\end{array}\right., $$
is a solution of the TMTMP.

\begin{center}
{\large\bf The Nevanlinna-type formula for the truncated matrix trigonometric moment problem.}
\end{center}
\begin{center}
{\bf S.M. Zagorodnyuk}
\end{center}

This paper is a continuation of our previous investigation on the truncated matrix trigonometric moment problem
in Ukrainian Math. J., 2011, 63, no.6, 786-797. In the present paper we obtain a Nevanlinna-type formula
for this moment problem in a general case. We only assume that we have more than one moment, the moment problem is
solvable and the problem has more than one solution.
The coefficients of the corresponding matrix linear  fractional transformation are
explicitly expressed by the prescribed moments.
Easy conditions for the determinacy
of the moment problem are given.

}

\begin{thebibliography}{1}

\bibitem{cit_500_Z}
Zagorodnyuk S.~M. The truncated matrix trigonometric moment problem: the operator approach //
Ukrainian Math. J. - 2011. - \textbf{63}, no. 6. - P. 786-797.

\bibitem{cit_4500_A}
Ando T., Truncated moment problems for operators //
Acta Scientarum Math., (Szeged).- 1970.- 31, no. 4.- P.319-334.

\bibitem{cit_6000_M}
Berezanskii Ju. M. Expansions in Eigenfunctions of Selfadjoint Operators. - Amer. Math. Soc., Providence, RI, 1968.
(Russian edition: Naukova Dumka, Kiev, 1965).

\bibitem{cit_7000_G}
Gantmacher F.  R. Theory of Matrices. - Nauka, Moscow, 1967 (in Russian).



\end{thebibliography}
\end{document}